\documentclass[11pt,a4paper]{amsart}

\usepackage[T1]{fontenc}
\usepackage[utf8]{inputenc}
\usepackage[english]{babel}
\usepackage{lmodern}
\usepackage[a4paper,margin=27mm]{geometry}
\usepackage{amsmath,amssymb,amsthm,mathtools}
\usepackage{microtype}
\usepackage[hidelinks]{hyperref}

\newcommand{\N}{\mathbb N}
\newcommand{\Cantor}{2^\omega}
\newcommand{\Baire}{\omega^\omega}
\newcommand{\GEN}{\mathsf{GEN}}
\newcommand{\OneGEN}{1\text{-}\mathsf{GEN}}
\newcommand{\FIP}{\mathsf{FIP}}
\newcommand{\ran}{\operatorname{ran}}
\newcommand{\Split}{\operatorname{Split}}
\newcommand{\restr}{\upharpoonright}
\newcommand{\W}{\mathrm W}

\theoremstyle{plain}
\newtheorem{theorem}{Theorem}
\newtheorem{lemma}[theorem]{Lemma}
\newtheorem{corollary}[theorem]{Corollary}
\theoremstyle{remark}
\newtheorem*{remark}{Remark}

\title[1-Genericity is not Weihrauch reducible to FIP]{
1-Genericity Is Not Weihrauch Reducible\\
to the Finite Intersection Principle}
\author{Patrizio Cintioli}
\address{Mathematics Division, School of Science and Technology,
University of Camerino, 62032 Camerino, Italy}
\email{patrizio.cintioli@unicam.it}

\begin{document}

\begin{abstract}
This paper shows that the problem of producing a Cohen $1$-generic real
relative to a given oracle is not Weihrauch reducible to the finite
intersection principle (FIP). More strongly, no nontrivial computable
FIP instance admits a single Turing functional that produces an
absolute $1$-generic real from every enumeration of every maximal
centered set of indices. The argument combines a $1$-generic real
that computes no perfect tree all of whose paths are $1$-generic
with local constancy and a parametrization of maximal centered sets.
This answers Question~9.10.25 of Dzhafarov and Mummert.
\end{abstract}
\maketitle

\section{Introduction}

The \emph{finite intersection principle} (FIP) asks, for every
nontrivial countable family of subsets of $\N$, for an
inclusion-maximal subfamily in which every finite collection of
at least two distinct members has nonempty intersection. Its relationship with
Cohen genericity has been studied in both computability theory and
reverse mathematics. Cholak, Downey, and Igusa~\cite{CDI} proved
that a Turing degree is FIP-bounding if and only if it computes
a Cohen $1$-generic; they also established the equivalence of FIP
and $1$-GEN over $\mathsf{RCA}_0$. Immediately after Theorem~1.6,
they contrasted the uniform extraction
available in the earlier $\Delta^0_2$ setting with
their unrestricted argument, which uses infinitely
many candidate Turing functionals. They also
reported preliminary evidence that no single
uniform extraction procedure could work in the
general case. The present result makes this
nonuniformity obstruction precise in the setting
of Weihrauch reducibility.

The uniform computational content of $1$-genericity
was investigated by Brattka, Hendtlass, and
Kreuzer~\cite{BHK}. The Weihrauch-theoretic study
of FIP also appears in the work of Dzhafarov,
Solomon, and Yokoyama~\cite{DSY}, concerning the
first-order parts of computational problems.
The question here is whether an entire generic
real can be extracted uniformly.
Dzhafarov and Mummert~\cite[Question~9.10.25]{DMbook}
explicitly asked whether $1$-genericity is Weihrauch
reducible to FIP. The present paper gives a negative
answer. The obstruction has a stronger form: a single
postprocessor cannot succeed on every solution
presentation of \emph{any fixed} nontrivial computable
FIP instance.

Let $\GEN$ be the class of absolute Cohen $1$-generic
reals in $\Cantor$.
Recall that \(G\in\Cantor\) is \(1\)-generic if, for every
c.e.\ set \(U\subseteq 2^{<\omega}\), at least one of the
following alternatives holds:
\[
(\exists n)\,G\restr n\in U,
\qquad
(\exists n)(\forall\tau\succeq G\restr n)\,\tau\notin U.
\]

A nontrivial computable family is a sequence
\(\mathcal A=(A_i)_{i\in\N}\) of subsets of \(\N\),
with the relation \(x\in A_i\) uniformly computable in \(i,x\),
and with at least one nonempty \(A_i\).

\begin{theorem}\label{thm:main}
For no nontrivial computable family \(\mathcal A\)
does there exist a Turing functional \(\Phi\) such that, for every
maximal centered set of indices \(M\) and every
enumeration \(J\) of all of \(M\), the real \(\Phi^J\)
is total and belongs to \(\GEN\).
\end{theorem}

An enumeration is a total sequence, possibly with repetitions,
whose range is exactly the set being enumerated. The required
functional must work for all orders and repetitions: no
permutation-invariance of its output is assumed.

As an immediate consequence of Theorem~\ref{thm:main}, proved in
Section~\ref{sec:weihrauch}, we obtain
\[
  \OneGEN\nleq_{\W}\FIP,
\]
where $\OneGEN$ is the multivalued problem that assigns to an oracle
$X$ the collection of reals that are $1$-generic relative to $X$.
This distinction between reverse-mathematical implication and
Weihrauch reducibility is essential: the latter requires a single
uniform pair of preprocessing and postprocessing procedures.

The proof of the stronger theorem makes essential use of the
following result.

\begin{theorem}[Hirschfeldt--Jockusch--Schupp,
Theorem 7.23]\label{thm:hjs}
There exists a \(1\)-generic real \(G\) that computes no
perfect tree all of whose paths are \(1\)-generic.
\end{theorem}

This is \cite[Theorem~7.23]{HJS}; its proof uses Kumabe's
construction of a $1$-generic degree with a strong minimal
cover~\cite{Kumabe}. The perfect trees below are represented by
splitting maps as in \cite[Definition~3.20]{HJS}, with computability
relative to the designated oracle.

The main idea is to parametrize centered sets by reals so that
$1$-generic parameters yield maximal centered sets and any finite
centered extension can be realized above any parameter prefix.
Theorem~\ref{thm:hjs} then provides a parameter for which our
hypothetical postprocessor is locally constant on one cone of
enumerations. Genericity promotes this local constancy to a
neighborhood of the parameter, and the extension property
transfers it to every maximal centered set extending a fixed
finite condition. A computable-output argument finally yields
a computable $1$-generic, a contradiction.

\section{Finite conditions and representation of solutions}

The relevant finite conditions are finite sets of indices whose
corresponding family members have nonempty intersection. A solution
is represented by an enumeration of a maximal collection of such
compatible indices; the order and repetitions are unrestricted.
The formulation in
\cite[Definition~9.10.17]{DMbook}
requires nonempty intersections for all nonempty
finite sets of indices, including singletons.
The centered sets used here satisfy this stronger
condition. Consequently, the solutions produced
by the construction are valid under either convention.

Fix \(\mathcal A\), and set
\[
\mathcal C=
\left\{
F\subseteq_{\mathrm{fin}}\N:
\bigcap_{i\in F}A_i\ne\varnothing
\right\},
\]
with the convention that \(\varnothing\in\mathcal C\).
The class \(\mathcal C\) is c.e.\ and downward closed.
The requirement on singletons makes every selected family member
nonempty; this stronger finite condition is convenient for the proof
and still produces solutions to the usual FIP problem.
For \(F\ne\varnothing\), a positive witness to
\(F\in\mathcal C\) is a number \(w\) belonging to every
\(A_i\) with \(i\in F\).

A set \(M\subseteq\N\) is \emph{centered} if each of its
finite subsets belongs to \(\mathcal C\).
It is \emph{maximal} if it is inclusion-maximal among
centered sets.

Every centered set extends to a maximal one.
It suffices to consider all indices in increasing order and add
each index whenever doing so preserves centeredness.
This is an existence construction, not a computable
algorithm.

For a maximal set \(M\), define
\[
E_M=\{J\in\Baire:\ran(J)=M\}.
\]
A maximal centered set is nonempty because the family is
nontrivial. Every finite string with values in \(M\) can be
extended to an element of \(E_M\), even if \(M\) is finite,
since repetitions are allowed.

For a finite string \(\sigma\), let \([\sigma]\)
denote the cylinder of sequences extending it, in Baire
or Cantor space according to the context.
When we write that \(\Phi^\sigma\) determines a binary
string \(u\), we mean a finite computation of the first
\(|u|\) bits that uses only positions available in \(\sigma\).

We write
\[
\Split_\Phi(\tau_0,\tau_1)
\]
for the property that there exist incompatible binary strings
\(u_0,u_1\) such that \(\Phi^{\tau_i}\) determines \(u_i\),
for \(i=0,1\). This property is c.e.

\section{Local constancy and computable output}

The first lemma converts constancy along the enumerations of each
maximal extension of one finite condition into a computable output.
The construction does not require a computable maximal centered set.

\begin{lemma}\label{lem:constant}
Suppose that \(\Phi\) is total, with values in \(\Cantor\),
on all enumerations of all maximal centered sets of
\(\mathcal A\).
Let \(c\in\mathcal C\) be nonempty, and let
\(\sigma\in\omega^{<\omega}\) satisfy \(\ran(\sigma)\subseteq c\).
If, for every maximal \(M\supseteq c\), the restriction
of \(\Phi\) to \(E_M\cap[\sigma]\) is constant, then
some output of \(\Phi\) on a maximal solution
is computable.
The constant value may depend on \(M\).
\end{lemma}

\begin{proof}
Let \(q_0\) be the concatenation of \(\sigma\) with the
increasing list of the elements of \(c\); thus \(\ran(q_0)=c\).

Computably construct increasing strings \(q_n\), with
centered ranges, such that \(|q_n|\ge n\) and
\(\Phi^{q_n}\) determines at least \(n\) bits.

Given \(q_n\), search for \(q_{n+1}\succeq q_n\) with
centered range, length at least \(n+1\), and a convergent
computation of at least \(n+1\) bits.
The search is c.e., using a positive witness to
\(\ran(q_{n+1})\in\mathcal C\).

The search terminates: \(\ran(q_n)\) extends to a maximal
set \(N\), and an enumeration of \(N\) beginning with \(q_n\)
makes \(\Phi\) total. A sufficiently long initial segment
of that enumeration supplies the required witness.

The computations are compatible because the input strings
extend one another. The limit of their outputs is a
computable real \(g\).

The set
\[
S=\bigcup_n\ran(q_n)
\]
is centered. Fix a maximal \(M\supseteq S\).
For every \(n\), there exists \(J_n\in E_M\) extending \(q_n\).
All these \(J_n\) extend \(\sigma\), so their outputs
\(\Phi^{J_n}\) coincide by hypothesis.
The common value extends \(g\restr n\) for every \(n\),
and hence equals \(g\).
\end{proof}

\begin{remark}
We do not construct a computable maximal solution.
We construct a c.e.\ centered set, extend it abstractly
to a maximal set, and use constancy on the cone
of enumerations of that maximal set.
The finite constants \(c,\sigma\) can be hard-coded
into the program computing \(g\).
\end{remark}

\section{From output variation to a perfect tree}

If the outputs vary above every finite enumeration stem, this variation
can be organized into a perfect binary tree of outputs. Positive access
to the given enumeration suffices to compute the tree.

\begin{lemma}\label{lem:perfect}
Let \(M\ne\varnothing\), and suppose that \(\Phi\) is total
on \(E_M\), with values in a class \(\mathcal Q\subseteq\Cantor\).
At least one of the following alternatives holds:
\begin{enumerate}
\item there exists a string \(\sigma\) with values in \(M\)
such that \(\Phi\) is constant on \(E_M\cap[\sigma]\);
\item every \(J\in E_M\) computes a perfect tree
\(T\subseteq2^{<\omega}\) with \([T]\subseteq\mathcal Q\).
\end{enumerate}
\end{lemma}

\begin{proof}
Suppose that the first alternative fails, and fix
\(J\in E_M\).
Above every finite string \(p\) with values in \(M\),
there are two enumerations of \(M\) beginning with \(p\)
and having distinct outputs.
Hence there exist two finite extensions \(p_0,p_1\) of \(p\),
still with values in \(M\), such that
\(\Split_\Phi(p_0,p_1)\).

This pair can be searched for effectively using \(J\):
membership of the finitely many values of \(p_0,p_1\) in \(M\)
can be verified positively using \(J\), while the splitting
is witnessed by finite computations.
The characteristic function of \(M\) is not used.

Using \(J\) as an oracle, construct input strings \(p_\alpha\)
and output strings \(u_\alpha\), for
\(\alpha\in2^{<\omega}\), with the following properties:
\begin{enumerate}
\item if \(\alpha\preceq\beta\), then
\(p_\alpha\preceq p_\beta\) and \(u_\alpha\preceq u_\beta\);
\item \(u_{\alpha^\frown0}\) and \(u_{\alpha^\frown1}\)
are incompatible and properly extend \(u_\alpha\);
\item \(\Phi^{p_\alpha}\) determines \(u_\alpha\), and
\(|u_\alpha|\ge|\alpha|\);
\item in passing from level \(n\) to level \(n+1\),
before searching for a splitting, append \(J(n)\)
to every input stem.
\end{enumerate}

Start with \(p_\varnothing=u_\varnothing=\varnothing\).
The search described above guarantees each step.
The children's outputs can be extended as far as needed
using totality on completions in \(E_M\).

For every \(B\in\Cantor\), the limit
\[
f_B=\bigcup_n p_{B\restr n}
\]
is an enumeration of all of \(M\): it uses only elements
of \(M\), includes \(J(n)\) for every \(n\), and has
infinite length. Therefore
\[
\bigcup_n u_{B\restr n}=\Phi^{f_B}\in\mathcal Q.
\]

The prefix closure of the strings \(u_\alpha\)
is a perfect tree \(T\).
It is computable in \(J\), not merely c.e.\ in \(J\):
to decide whether a string \(v\) of length \(m\)
belongs to \(T\), it suffices to construct level \(m\)
of the outputs and check whether \(v\) is a prefix of one
of them. All these outputs have length at least \(m\).

To see that every path through $T$ is one of the
outputs described above, put
\[
L_n=\max\{|u_\alpha|:|\alpha|=n\}.
\]
Every path through $T$ extends exactly one of the
strings $u_\alpha$ at level $n$. These choices are
coherent across levels and therefore determine a
unique branch $B\in2^\omega$. The path is consequently
equal to $\bigcup_n u_{B\restr n}=\Phi^{f_B}$.
Thus $[T]\subseteq\mathcal Q$.

Moreover, the map $S(\alpha)=u_\alpha$ directly
satisfies the representation of perfect trees used
in \cite[Definition~3.20]{HJS}: the children are
incompatible and properly extend the parent.
Thus no noncomputable conversion between
representations is needed.
\end{proof}

\section{A map from generics to maximal sets}

Forcing with strings having centered ranges, and its connection
to genericity, appear in
\cite[Section 2, Lemma 2.1]{DDGT}.
Here we make explicit a total parametrization
on \(\Cantor\) with the required local universality.

\begin{lemma}\label{lem:universal}
There exist a total functional \(R\mapsto J_R\), from
\(\Cantor\) to \(\Baire\), and a computable function
\(\xi\mapsto c_\xi\in\mathcal C\), defined on finite binary
strings, with the following properties:
\begin{enumerate}
\item \(c_\xi\ne\varnothing\), and
\(c_\xi\subseteq c_\eta\) whenever \(\xi\preceq\eta\);
\item \(M_R:=\ran(J_R)=\bigcup_n c_{R\restr n}\)
is centered;
\item if \(R\) is \(1\)-generic, \(M_R\) is maximal;
\item for every \(\xi\) and every \(D\in\mathcal C\)
with \(D\supseteq c_\xi\), there exists \(\eta\succeq\xi\)
with \(c_\eta=D\).
In particular, \(D\subseteq M_R\) for every \(R\succeq\eta\).
\end{enumerate}
\end{lemma}

\begin{proof}
Find \(i_0\) with \(A_{i_0}\ne\varnothing\) by searching
for a pair \((i_0,w)\) with \(w\in A_{i_0}\).
Set \(F_0=\{i_0\}\).

Use a computable bijective coding of triples \((D,w,r)\),
where \(D\) is a finite set and \(r\) is a free parameter.
If \(k\) codes \((D,w,r)\), define
\[
e(F,k)=
\begin{cases}
D,&F\subseteq D\ \text{and}\ (\forall i\in D)\,w\in A_i,\\
F,&\text{otherwise}.
\end{cases}
\]
The function \(e\) is total and computable.
If \(F\) is centered, \(e(F,k)\) is a centered extension
of \(F\).
Every \(D\in\mathcal C\) with \(F\subseteq D\) occurs
as \(e(F,k)\) for arbitrarily large values of \(k\),
by choosing a witness \(w\) and varying \(r\).

Read the completed blocks of \(R\) in the form
\[
0^{k_0}1\,0^{k_1}1\,0^{k_2}1\cdots
\]
and update \(F_{n+1}=e(F_n,k_n)\).
For a finite string \(\xi\), let \(c_\xi\)
be the last \(F_n\) obtained from its completed blocks.
Ignore any incomplete final block;
in particular \(c_\varnothing=F_0\).

To make \(J_R\) total for every \(R\), concatenate,
for \(n=0,1,\ldots\), the increasing, finite, nonempty
lists of elements of \(c_{R\restr n}\).
This enumerates their union computably in \(R\),
even if \(R\) contains only finitely many \(1\)'s.
The first two properties follow from monotonicity.

We verify local universality.
If \(\xi\) ends with an incomplete block of \(a\) zeros,
and \(D\supseteq c_\xi\) is centered, choose \(k\ge a\)
with \(e(c_\xi,k)=D\).
Then
\[
\eta=\xi^\frown 0^{k-a}1
\]
satisfies \(c_\eta=D\).
The free parameter ensures that every prefix,
even one ending in the middle of a block, allows every
future centered extension.

For maximality, fix \(i\) and set
\[
U_i=\{\xi:i\in c_\xi\}.
\]
The set \(U_i\) is computable and upward closed under
extension. If \(R\) is \(1\)-generic, either it meets \(U_i\),
in which case \(i\in M_R\), or there is a \(\xi\prec R\)
no extension of which belongs to \(U_i\).

In the second case \(c_\xi\cup\{i\}\notin\mathcal C\):
otherwise, local universality would yield an extension
of \(\xi\) in \(U_i\).
Thus \(i\) cannot be added to \(M_R\).
This proves maximality.
\end{proof}

\section{Proof of the main theorem}

\begin{proof}[Proof of Theorem~\ref{thm:main}]
Suppose, toward a contradiction, that a nontrivial
computable family \(\mathcal A\) and a functional \(\Phi\)
satisfy the property ruled out by the statement.

Choose \(G\) as in Theorem~\ref{thm:hjs}.
Apply Lemma~\ref{lem:universal}, and put
\[
J=J_G,\qquad M=M_G.
\]
Then \(J\le_T G\), and \(M\) is maximal.

Apply Lemma~\ref{lem:perfect} to \(M,\Phi\)
and \(\mathcal Q=\GEN\).
The second alternative is impossible:
it would give a perfect tree of \(1\)-generic reals
computable in \(J\), and hence in \(G\).
Thus there exists a finite string \(\sigma\),
with \(\ran(\sigma)\subseteq M\), such that \(\Phi\)
is constant on \(E_M\cap[\sigma]\).

For this \(\sigma\), define
\[
F_\sigma=
\left\{
R\in\Cantor:
\begin{array}{l}
\text{there do not exist }\tau_0,\tau_1\succeq\sigma
\text{ such that}\\
\ran(\tau_0)\cup\ran(\tau_1)\subseteq M_R
\text{ and }\Split_\Phi(\tau_0,\tau_1)
\end{array}
\right\}.
\]
To make its complexity explicit, put
\[
V_\sigma=
\left\{
\zeta\in2^{<\omega}:
\begin{array}{l}
\text{there exist }\tau_0,\tau_1\succeq\sigma
\text{ such that}\\
\ran(\tau_0)\cup\ran(\tau_1)\subseteq c_\zeta
\text{ and }\Split_\Phi(\tau_0,\tau_1)
\end{array}
\right\}.
\]

The set \(V_\sigma\) is c.e.\ and upward closed under
extension: \(c_\zeta\) is finite and computable,
and a splitting has finite computational witnesses.
By monotonicity of the \(c_\zeta\),
\[
R\notin F_\sigma
\iff
(\exists n)\,R\restr n\in V_\sigma.
\]
Indeed, finitely many elements of
\(M_R=\bigcup_n c_{R\restr n}\) all occur
in the same \(c_{R\restr n}\) for some \(n\).
Therefore
\[
F_\sigma=\Cantor\setminus
\bigcup_{\zeta\in V_\sigma}[\zeta]
\]
is a \(\Pi^0_1\) class.
The definition applies to every \(R\), even when \(M_R\)
is not maximal, and does not assume totality of
\(\Phi^{J_R}\) on those parameters.

Constancy on \(E_{M_G}\cap[\sigma]\) implies
\(G\in F_\sigma\).
Otherwise, the two finite prefixes with incompatible
outputs could be extended to two complete enumerations
of \(M_G\) belonging to the same cone.

A \(1\)-generic real belonging to a \(\Pi^0_1\) class
belongs to its interior.
Applying genericity to \(V_\sigma\), find
\(\xi\prec G\) with
\[
[\xi]\subseteq F_\sigma.
\]
Extend \(\xi\) further along \(G\) until
\(\ran(\sigma)\subseteq c_\xi\), and put \(c=c_\xi\).
This is possible because
\(M_G=\bigcup_n c_{G\restr n}\) and \(\ran(\sigma)\)
is finite.

We claim that, for every maximal \(N\supseteq c\),
\(\Phi\) is constant on \(E_N\cap[\sigma]\).
If not, two enumerations in that cone would have
distinct outputs.
There would then be finite prefixes
\(\tau_0,\tau_1\succeq\sigma\), with values in \(N\),
such that \(\Split_\Phi(\tau_0,\tau_1)\).

The set
\[
D=c\cup\ran(\tau_0)\cup\ran(\tau_1)
\]
is centered, as a finite subset of the same
centered set \(N\).
By local universality there exists \(\eta\succeq\xi\)
with \(c_\eta=D\).
For every \(R\succeq\eta\), the same prefixes
witness \(R\notin F_\sigma\), contradicting
\[
[\eta]\subseteq[\xi]\subseteq F_\sigma.
\]
This proves the claim.

Lemma~\ref{lem:constant}, applied to \(c,\sigma\),
now provides a computable output \(g\) of \(\Phi\)
on a maximal solution.
By hypothesis, \(g\) should be \(1\)-generic.

No computable real is \(1\)-generic:
if \(g\) is computable, the computable set
\[
U=\{(g\restr n)^\frown(1-g(n)):n\in\N\}
\]
contains no prefix of \(g\), but contains an extension
of every prefix of \(g\).
This is the required contradiction.
\end{proof}

\section{Consequence for Weihrauch reducibility}\label{sec:weihrauch}

The preceding theorem already fixes one computable instance and
rules out all uniform postprocessors. This immediately excludes
an ordinary Weihrauch reduction by specializing its input to a
computable oracle.

Regard \(\OneGEN\) as the problem which, given
an oracle \(X\), asks for a real that is \(1\)-generic
relative to \(X\).
FIP associates with every nontrivial family a
maximal subfamily having the finite intersection
property, represented by a sequence in which
repetitions are permitted.

\begin{corollary}
\[
\OneGEN\nleq_{\W}\FIP.
\]
\end{corollary}

\begin{proof}
Suppose a pair of computable functionals
\(\Delta,\Psi\) witnesses
\(\OneGEN\leq_{\W}\FIP\).
Fix the computable oracle \(X=\varnothing\)
as the initial instance.
The preprocessor produces a nontrivial
computable family
\[
\mathcal A=\Delta(\varnothing).
\]

For every maximal centered set of indices \(M\)
and every \(J\in E_M\), the sequence
\[
(A_{J(n)})_{n\in\N}
\]
is an FIP solution computable from \(J\).
Composing this transformation with the postprocessor
\(\Psi\) and the computable parameter \(X=\varnothing\),
we obtain a functional \(\Phi\) that produces an
absolute \(1\)-generic from all these enumerations.
This contradicts Theorem~\ref{thm:main}.
\end{proof}

\section{Representations and nonuniformity}

The separation is for \emph{ordinary} Weihrauch reducibility,
not merely its strong variant: fixing the computable initial
instance in the preceding proof absorbs the original input
accessible to the postprocessor.

Indices and sets need not be in one-to-one correspondence:
repetitions among the $A_i$ are permitted. Maximality of a
centered set of indices gives maximality of the corresponding
family of sets, since any additional set that could be adjoined
would have an index that could also be adjoined. This remains true
when FIP is formulated using collections of at least two distinct
members: the selected sets are nonempty, so such a maximal centered
family is a valid FIP solution. The converse is not needed. We use only
the computable conversion
\[
  J\longmapsto (A_{J(n)})_{n\in\omega},
\]
not its converse. In particular, no decision procedure for
equality of computable sets is required; cf.\ the discussion of
solution presentations in~\cite[Section~1.1]{DDGT}.

The splitting argument uses the union of finite input ranges
only when those ranges are contained in a \emph{common}
centered set. The finite constants chosen in the proof,
including $\sigma$ and $c_\xi$, need not be found uniformly
from $\mathcal A$; once chosen, they can be hard-coded in the
computation of the final real. Thus the argument does not
conflict with the nonuniform computability-theoretic equivalence
proved in~\cite{CDI}. Their earlier observation concerning the
failure of uniformity is an important antecedent of the present
separation.

Dzhafarov and Mummert also asked in
\cite[Question~9.10.26]{DMbook} whether
$\mathsf{RCA}_0\vdash
2\mathrm{IP}\to1\text{-}\mathrm{GEN}$.
This question is answered affirmatively in a companion
manuscript~\cite{CintioliIP}.
There is no conflict between this reverse-mathematical
result and the Weihrauch nonreducibility established
here: provability of an implication over
$\mathsf{RCA}_0$ does not, by itself, yield a
Weihrauch reduction. In particular, the companion
proof uses a nonuniform distinction between an
already generic path and a case in which the jump
of the original instance can be decoded.

\section*{Acknowledgments}
The main result was discovered by ChatGPT 6 Pro (OpenAI),
which was also used in developing the proof and preparing
the manuscript. The manuscript was subsequently subjected
to separate AI-assisted mathematical audits using
Astra 6 Ultra (OpenAI) and Gemini 3.1 Deep Think
(Google DeepMind). The author takes full responsibility
for the correctness of the mathematical statements and
proofs, for the accuracy of the citations and attributions,
and for the final content of the manuscript.

\end{document}